\documentclass{article}
\usepackage{latexsym}
\usepackage{amssymb}
\usepackage{amsmath}
\usepackage{color}
\usepackage{amsfonts,amssymb}
\usepackage{mathrsfs}
\usepackage{dsfont}

\begin{document}
\title{\bf ${\bf -1}$  Krall-Jacobi polynomials}
\author{ Luc Vinet$^{\dag}$, Guo-Fu Yu$^{\ddag,\dag}$ and Alexei Zhedanov$^{\S}$\\
\\ $^{\dag}$
Centre de Recherches Math\'{e}matiques, Universit\'{e} de
Montr\'{e}al,
\\ C.P.6128, Centre-ville Station, Montr\'{e}al, Qu\'{e}bec, H3C 3J7, Canada\\
$^{\ddag}$Department of Mathematics, Shanghai Jiao Tong University, \\
Shanghai 200240, P.R.\ China \\
$^{\S}$ Donetsk Institute for Physics and Technology. Donetsk 83114,
Ukraine}
\date{}
\maketitle

\begin{abstract}
We study a family of orthogonal polynomials which satisfy (apart
from a 3-term recurrence relation) an eigenvalue equation involving
a third order differential operator of Dunkl-type. The orthogonality measure of these polynomials consists in the continuous measure of the little -1 Jacobi polynomials 
to which is added an arbitrary mass located at the point $x=0$, the middle of the orthogonality interval. This provides the first nontrivial example of Krall-type polynomials 
with a point mass inside the orthogonality interval. These polynomials
can be obtained by a Geronimus transform of the little
$q$-Jacobi polynomials in the limit $q=-1$.

\end{abstract}

\section{Introduction}
Significant advances have been realized in the characterization of
recently discovered families of $-1$ orthogonal polynomials (OPs) that can be obtained as $q \to -1$ limits of $q$ polynomials of the Askey scheme.
The striking feature of these OPs is that they are classical or
bispectral and  that they satisfy eigenvalue equations involving
Dunkl-type operators \cite{dunkl} in addition to the mandatory
3-term recurrence relation. They have arisen already in a number of
physical problems \cite{vz1}-\cite{gvz} and are connected with
Jordan algebras \cite{tvz2}.

At the top of the emerging $-1$ scheme are the Bannai-Ito
polynomials and their kernel partners (\cite{Chi},Ch.1, Sect. 7), the complementary Bannai-Ito
polynomials \cite{tvz3}. Both sets depend on $4$ parameters. The
Bannai-Ito polynomials are the eigensolutions of the most general
operator which is of first order in Dunkl shifts, (i.e. first order
in the operators $T$ and $R$ defined by $Tf(x)=f(x+1), Rf(x)=f(-x)$
on a function $f(x)$) and which stabilizes polynomials of given
degrees. The Bannai-Ito polynomials are positive definite and
orthogonal on a finite set of $N+1$ points. In the limits where
$N\rightarrow \infty$, they tend to the big $-1$ Jacobi polynomials
\cite{vz3} which are orthogonal on $[-1,-c]\bigcup [c,1]$. When
$c=0$, the little $-1$ Jacobi polynomials \cite{mop} arise as a
special case. A Bochner-type theorem establishes \cite{vz4} that the
big and little $-1$ polynomials are the only families of OPs
satisfying a differential-difference eigenvalue equation which is of
first order in Dunkl-type operators.

The Bannai-Ito polynomials and their descendants  all possess the
Leonard duality property. (In fact, this led to their initial
identification \cite{Ba}.) The dual $-1$ Hahn polynomials
\cite{tvz4} together with the generalized Gegenbauer \cite{bel} and
Hermite polynomials \cite{rosenblum, BenGaied} are also bispectral
but, obeying eigenvalue equations of second order in Dunkl
operators, they fall beyond the scope of Leonard duality.

We continue here the exploration of $-1$ orthogonal polynomials in
that vein and look for another class of $-1$ OPs verifying a higher
differential-difference equation.

In the wake of Krall's classification of OPs satisfying a fourth
order differential equations \cite{Krall}, it is appreciated that
the addition of discrete masses to the measure leads to OPs
verifying higher order equations \cite{koor,koekoek, kwon}.

In this connection, we present here a generalization of the little
$-1$ Jacobi polynomials with the following features: these
orthogonal polynomials obey a differential-difference equation of
third order in Dunkl operators and a mass is located at the middle
of the orthogonality interval.

The outline of the paper is the following. In section 2, we offer a
brief review of useful results on little q-Jacobi polynomials. In
section 3, we introduce the generalized little q-Jacobi polynomials
\cite{GqJ} that are eigensolutions of higher order q-difference
equations. In section 4, focusing for definiteness on one of the
simpler cases, we obtain and characterize a set of generalized
little $-1$ Jacobi polynomials by taking an appropriate q
$\rightarrow -1$ limit of certain polynomials of the preceding
section. The paper ends with concluding remarks.

\section{Little $q$-Jacobi polynomials}
\setcounter{equation}{0} The monic little q-Jacobi polynomials are
defined as
\begin{eqnarray}
P_n(x;a,b)=(-1)^n\frac{q^{n(n-1)/2}(aq;q)_n}{(abq^{n+1};q)_n}
{}_2\phi_1\left(\begin{matrix}q^{-n},abq^{n+1}\\ aq
\end{matrix}| q x \right)
\end{eqnarray}
where $(a;q)_n=(1-a)(1-aq)\cdots (1-aq^{n-1})$ is the q-shifted
factorial and ${}_2\phi_1$ denotes the q-hypergeometric function.

The orthogonality relation is
\begin{eqnarray}
\sum_{k=0}^\infty w_k P_{n}(q^k;a,b)P_m(q^k;a,b)=h_n\delta_{nm},
\end{eqnarray}
where $h_n$ are appropriate normalization constants, and the
normalized weight function is
\begin{eqnarray}
w_k=\frac{(aq;q)_{\infty}}{(abq^2;q)_{\infty}}\frac{(bq;q)_k(aq)^k}{(q;q)_k}.
\end{eqnarray}
It is assumed that $|q|<1, \: |aq|<1, \: |bq|<1$.  The expansion coefficients
of the little q-Jacobi polynomials in
$P_n(x;a,b)=\sum\limits_{s=0}^n B_n^{(s)}x^{n-s}$ are
\begin{eqnarray}
B_n^{(s)}=b^{-s}\frac{(q^{-n},a^{-1}q^{-n};q)_s}{(q,a^{-1}b^{-1}q^{-2n};q)_s},
\end{eqnarray}
where $(a_1,a_2, \dots, a_k;q)_n = (a_1,q)_n (a_2,q)_n \dots (a_k,q)_n$.

It is known that the little q-Jacobi polynomials satisfy a
second-order difference equation \cite{KS}.

Introduce the functions of the second kind,
\begin{eqnarray}
Q_n(z)=\int_a^b \frac{P_n(x)w(x)}{z-x}dx,\label{qz}
\end{eqnarray}
where $w(x)$ is assumed to be the normalized weight function, i.e.
\[
\int_a^b w(x)dx=1.
\]
The values of the functions $Q_n(z)$, at $z=0$ (an accumulation
point of the orthogonality measure), are :
\begin{eqnarray}
Q_n(0;a,b)=-\sum_{k=0}^\infty \frac{P_n(q^k;a,b)w_k}{q^k}.
\end{eqnarray}
Using the $q-$ binomial theorem and the $q-$Saalsch\"{u}tz formula
(see, e.g.,\cite{Gasper}) we have
\begin{eqnarray}
Q_n(0;a,b)=(-1)^{n+1}a^nq^{n(n-1)/2}\frac{1-abq}{1-a}
\frac{(q;q)_n(bq;q)_n}{(abq;q)_n(abq^{n+1};q)_n}.
\end{eqnarray}

\section{Transformed $q-$Jacobi polynomials}
\setcounter{equation}{0} Let $P_n(x)$ be orthogonal polynomials with
measure localized on the interval $[a,b]$. Let $w(x)$ be the
corresponding normalized unit weight function and $Q_n(z)$ be
defined by \eqref{qz}. Let finally $c$ be a point beyond the
orthogonality interval $[a,b]$ such that $Q_n(c)$ exists.

Consider the Geronimus transformation \cite{Ge1,Ge2} of the
polynomials $P_n(x)$ at the point $x=c$
\begin{eqnarray}
\tilde{P}_n(x)=\mathcal{G}(c){P_n(x)}=P_n(x)-\frac{\Phi_n}{\Phi_{n-1}}P_{n-1}(x),n=1,2,\cdots,
\tilde{P}_0(x)=1, \label{GT}
\end{eqnarray}
where
\begin{eqnarray}
\Phi_n=Q_n(c)+M P_n(c),
\end{eqnarray}
with $M$ some arbitrary real parameter.

We note that if $a=q^j, j=1,2,3,\cdots$, then
\begin{eqnarray}
&&\Phi_n=Q_n(0;q^j,b)+M P_n(0;q^j,b)\nonumber\\
&&\quad=(-1)^nq^{n(n-1)/2}\frac{(q^{j+1};q)_n}{(bq^{n+j+1};q)_n}
\left(
M-q^{nj}\frac{(1-bq^{j+1})(bq;q)_j(q;q)_j}{(1-q^j)(q^{n+1};q)_j(bq^{n+1};q)_j}
\right).\label{phiexp}
\end{eqnarray}

The weight function $\tilde{w}(x)$ of the polynomials
$\mathcal{G}(c)\{P_n(x)\}$ is
\begin{eqnarray}
\tilde{w}(x)=\kappa\left( \frac{w(x)}{x-c}-M \delta(x-c)\right),
\end{eqnarray}
where $\kappa$ is an appropriate normalization constant. The
Geronimus transformation thus inserts a concentrated mass at the
point $x=c$. The value of this mass depends on the parameter $M$.
Now take for $P_n(x)$ the little q-Jacobi polynomials with
$a=q^j,j=1,2,3,\cdots,$ and perform the Geronimus transformation
(\ref{GT}) with $\Phi_n$ given by (\ref{phiexp}). (In this case
$c=0$.)

The weight function $\tilde{w}(x)$ for the polynomials
$\mathcal{G}(0)\{P_n(x)\}$ is
\begin{eqnarray}
\tilde{w}(x)=\kappa \left( \sum_{k=0}^\infty
\tilde{w}_k\delta(x-q^k)-M\delta(x) \right),\label{wfun}
\end{eqnarray}
where
\begin{eqnarray}
\tilde{w}_k=\frac{(q^{j+1};q)_\infty}{(bq^{j+2};q)_\infty}\frac{(bq^k;q)_kq^{jk}}{(q;q)_k}.
\end{eqnarray}
The coefficients $B_n^{(s)}$ in the expansion
\begin{eqnarray}
\tilde{P}_n(x)=\mathcal{G}(0)\{P_n(x;q^j,b)\}=\sum_{s=0}^{n}B_n^{(s)}x^{n-s}
\label{TP}
\end{eqnarray}
have been given in \cite{GqJ}. In the same paper, it was shown that there exists an operator of the form  
\begin{eqnarray}
\mathcal{L}_q=\sum_{k=0}^{2N}a_k(q^{-N}x)T^{-N}\mathcal{D}_q^k,\label{Lq}
\end{eqnarray}
with
\begin{eqnarray}
a_k(x)=\sum_{s=0}^k\alpha_{ks}x^s,\quad k=0,1,\cdots,2N,\label{ak}
\end{eqnarray}
$T$  the $q-$shift operator and $\mathcal{D}_q$ the $q-$derivative
operator, such that the orthogonal polynomials $P_n(x)$ satisfy an
eigenvalue equation of the kind
\begin{eqnarray}
\mathcal{L}_qP_n(x)=\lambda_n P_n(x).\label{Lqe}
\end{eqnarray}

Consider the action of the operator $\mathcal{L}_q$ upon the
monomials $x^n$. From (\ref{Lq}) and (\ref{ak}) we get

\begin{eqnarray}
\mathcal{L}_q x^n=\sum_{s=0} A_n^{(s)} x^{n-s},\label{a3}
\end{eqnarray}
where
\begin{eqnarray}
A_n^{(s)}=q^{N(s-n)}[n][n-1]\cdots[n-s+1]\pi_s(q^n),
\end{eqnarray}
and
\begin{eqnarray}
\pi_s(q^n)=\alpha_{s0}+\sum_{i=1}^{2N-s}\alpha_{s+i,i}[n-s][n-s-1]\cdots[n-s-i+1]
\end{eqnarray}
are polynomials in $z=q^n$ of degree not exceeding $2N-s$. It is
clear that
\begin{eqnarray}
A_n^{(s)}=0,\quad s>2N
\end{eqnarray}
(see \cite{GqJ}). The coefficients $A_n^{(s)}$ completely
characterize the operator $\mathcal{L}_q$ and $A_n^{(s)}$ are called
the representation coefficients of the operator $\mathcal{L}_q$.

The coefficients $A_n^{(s)}$ for the $q$-difference operator
$\mathcal{L}_q$ that has the polynomials $\tilde{P}_n(x)$ (\ref{TP})
as eigenfunctions have been constructed in \cite{GqJ}, they are the
following
\begin{eqnarray}
&&A_n^{(0)}=\lambda_n=\frac{M
(q-1)q^{-n(j+1)-1}(q^n;q)_{j+1}(bq^n;q)_{j+1}}{1-q^{-j-1}}\nonumber\\
&&\qquad\qquad\qquad\qquad -(q^{-n}-1)(1-bq^{n+j})(b
q;q)_{j+1}(q;q)_{j-1},\label{a1} \\
&& A_n^{(1)}=(1-q^{-n})\left(M
q^{j(1-n)}(q^{n+1};q)_j(bq^n;q)_j(1-q^{n-1})-(q;q)_{j-1}(b
q;q)_{j+1}(1-q^{n+j-1}) \right),\\
&& A_n^{(s)}=M(q-1) q^{(s-n)(j+1)-1} \: \frac{(q^{-j};q)_{s-1} (q^{n-s};q)_{s+j+1}(bq^n;q)_{j-s+1}}{(q;q)_s} \quad \mbox{if} \quad s \ge 2 .\label{a2}
\end{eqnarray}
Note that from the above formula it follows $A_n^{(s)}=0$ if $s>j+1$

\section{Limit of the Krall-Jacobi polynomials as $q\rightarrow -1$}
\setcounter{equation}{0} In this section we construct the
$q\rightarrow -1$ limit of the coefficients $A_n^{(s)}$.  We take
$j=2$ and put
\begin{eqnarray}
q=-e^{\epsilon},\quad b=-e^{\beta \epsilon},\label{tr}
\end{eqnarray}
Note that $|qb|<1$ implies that $\beta+1>0$.

Substituting into the $A_n^{(s)}$ as given in (\ref{a1})-(\ref{a2})
and taking the limit $\epsilon \rightarrow 0$, we have
\begin{eqnarray}
&& \frac{A_n^{(0)}}{\epsilon^3}\rightarrow E_n^{(0)} = \left\{\begin{array}{ll}
-8M n(n+2)(n+1+\beta)+8n(\beta+1)(\beta+3) \quad & n \quad \texttt{even} \\
8M (n+1)(n+\beta)(n+2+\beta)-8(n+2+\beta)(\beta+1)(\beta+3) \quad &
n \quad \texttt{odd}\end{array}\right.\label{al1}
\\
&& \frac{A_n^{(1)}}{\epsilon^3}\rightarrow E_n^{(1)} = \left\{\begin{array}{ll}
8M n(n+2)(n+1+\beta)-8n(\beta+1)(\beta+3), \quad & n \quad
\texttt{even}
\\8(\beta+1)(\beta+3)(n+1)-8M(n^2-1)(n+\beta), \quad & n \quad
\texttt{odd}\end{array}\right.\\
&& \frac{A_n^{(2)}}{\epsilon^3}\rightarrow E_n^{(2)} = \left\{\begin{array}{ll}
8Mn(n+2)(n-2), \quad & n \quad \texttt{even}
\\-8M(n+1)(n-1)(n+\beta), \quad & n \quad
\texttt{odd}\end{array}\right.\\
&& \frac{A_n^{(3)}}{\epsilon^3}\rightarrow E_n^{(3)} = \left\{\begin{array}{ll}
-8Mn(n+2)(n-2), \quad & n \quad \texttt{even}
\\8M(n+1)(n-1)(n-3), \quad & n \quad
\texttt{odd} \quad.\end{array}\right.\label{al2}
\end{eqnarray}
and $E_n^{(s)}=0$ for $s=4,5,6,\dots$.

Consider the form of the $q$-difference equation (\ref{Lqe}) in this
limit. We divide both sides of (\ref{Lqe}) by $\epsilon^3$ and
introduce the operator $L_{\epsilon}$ which acts on the polynomials
$\tilde{P}_n(x)$ as
\begin{align}
L_{\epsilon}\tilde{P}_n(x)=\epsilon^{-3}\lambda_n \tilde{P}_n(x)
\end{align}
For monomials $x^n$, from (\ref{a3}) and (\ref{al1})-(\ref{al2}) we
have in the limit $\epsilon\rightarrow 0$,
\begin{align}
L_0 x^n=\lim\limits_{\epsilon\rightarrow 0}\frac{L_{\epsilon}}{\epsilon^3}x^n=& E_n^{(0)} x^n + E_n^{(1)} x^{n-1}  + E_n^{(2)} x^{n-2} + E_n^{(3)} x^{n-3}
\end{align}
This allows one to present the
operator $L_0=\lim\limits_{\epsilon\rightarrow
0}\frac{L_{\epsilon}}{\epsilon^3}$ in the form
\begin{align}
L_0=&(-8M+8Mx+8Mx^2-8M x^3)\partial_x^3R\nonumber\\
&+\left[-12M/x+24M+4\beta M+(36M+8\beta M)x-(12\beta
M+48M)x^2\right]\partial_x^2 R
\nonumber\\
&+(12Mx+4\beta Mx^2-12M/x-4\beta M)\partial_x^2+\left[(24M+16\beta
M-8\beta^2-32\beta-24)\right.\nonumber\\
&\left.+24M/x^2+(4\beta-12)M/x+(8\beta^2-36\beta M-48M-4\beta^2
M+32\beta+24)x\right]\partial_xR\nonumber\\
&+\left[(4\beta^2 M+12\beta M)x-(12+4\beta)M/x+24M+8\beta
M\right]\partial_x\nonumber\\
&+\left[12M/x^3+4\beta M/x^2+(12+4\beta^2+4\beta
M+16\beta)/x\right.\nonumber\\
&\left.+(8\beta M+4\beta^2
M-44\beta-24-4\beta^3-24\beta^2)\right](1-R),
\end{align}
where $R$ is the reflection operator $Rf(x)=f(-x)$. We thus have
that the polynomials $\tilde{P}_n^{(-1)}(x)$ are classical and
satisfy the eigenvalue equation
\begin{align}
L_0 \tilde{P}_n^{(-1)}(x)=\tilde{\lambda}_n \tilde{P}_n^{(-1)}(x),
\label{eigen}
\end{align}
where
\begin{eqnarray}
\tilde{\lambda}_n=E_n^{(0)}. 
\end{eqnarray}
The lower degree eigensolutions of (\ref{eigen}) can be obtained
directly as a check and as examples. We take
\begin{eqnarray}
\tilde{\lambda}_1=(\beta+1)(\beta+3)(16M-8(\beta+3)),
\end{eqnarray}
and find the first order polynomial solution
\begin{eqnarray}
\tilde{P}_1^{(-1)}(x)=x-1+\frac{2\beta-1-\alpha}{2\beta-3-\alpha}.
\end{eqnarray}
The second-order polynomial solution is
\begin{align}
&\tilde{\lambda}_2=(\beta+3)(16\beta+16-64M), \\
&\tilde{P}_2^{(-1)}(x)=x^2-\frac{2(4M-\beta-1)}{(5+\beta)(2M-\beta-1)}
x+\frac{2(\beta+1)}{(5+\beta)(2M-\beta-1)}.
\end{align}
and the third-order polynomial solution
\begin{align}
& \tilde{\lambda}_3=(3+\beta)(5+\beta)(32M-8-8\beta),\\
&
\tilde{P}_3^{(-1)}(x)=x^3-\frac{4(-2M+1+\beta)}{(7+\beta)(-4M+\beta+1)}x^2-\frac{4}{7+\beta}x
+\frac{8(1+\beta)}{(7+\beta)(5+\beta)(-4M+\beta+1)}.
\end{align}
Other polynomial eigensolutions can be obtained in the $q=-1$ limit
of \eqref{GT}. The weight function \eqref{wfun} has the following
moments
\begin{align}
\tilde{c}_n=k\left(\frac{(q^2;q)_n}{(bq^3;q)_n}-\frac{1-q^2}{1-bq^3}M\delta_{n,0}
\right)
\end{align}
with $k$ a different normalization than $\kappa$ in  \eqref{wfun}

Using the parametrization (\ref{tr}) and taking the limit $\epsilon
\rightarrow 0$, we can directly obtain from the above relation the
moments corresponding to the polynomials $\tilde{P}_n^{(-1)}$.
\begin{align}
\mu_{2n}=\mu_{2n-1}=k\left[\frac{(1)_n}{(\beta/2+3/2)_n}-\frac{2}{3+\beta}M\delta_{n,0}\right] = k \left[ \frac{n!}{(\beta/2+3/2)_n}+ M_0\delta_{n,0} \right],\qquad
n=1,2,\cdots,\label{nwet}
\end{align}
where $(x)_n=x(x+1)(x+2)\cdots (x+n-1)$ is the ordinary Pochhammer
symbol and $M_0 = -2M/(3+\beta)$. The parameter $M_0$ is the value of the concentrated mass inserted at the point $x=0$. It is seen that $M_0>0$ 
if $M<0$ (recall that condition $\beta>-1$ is assumed).

Indeed, when $M_0=0$ we have for the moments
\begin{align}
\mu_0=1, \; \mu_{2n}=\mu_{2n-1}=\frac{(1)_n}{(\beta/2+3/2)_n} , \quad n=1,2,\dots \label{mu_lJ} \end{align}
They correspond to a special case of the little -1 Jacobi polynomials \cite{mop} $P_n^{(\alpha,\beta)}(x)$ when $\alpha=1$.
Hence the moments (\ref{nwet}) are associated with the perturbation of the measure for the 
little -1 Jacobi polynomials $P_n^{(1,\beta)}(x)$ by the inserting of the concentrated mass $M_0$ at the point $x=0$.

It is then easily verified that
\begin{align}
w(x)=\tilde{k}\left(|x|(1-x^2)^{(\beta-1)/2}(1+x)
-\frac{4}{(1+\beta)(3+\beta)}M\delta(x)\right),
\end{align}
where
\begin{eqnarray}
\tilde{k}=\frac{\Gamma(\beta/2+3/2)}{\Gamma(\beta/2+1/2)}k=\frac{\beta+1}{2}k
\end{eqnarray}
is the orthogonality measure for these polynomials. Indeed we see
that
\begin{eqnarray}
\int_{-1}^{1}w(x)x^ndx=\mu_n,\qquad n=0,1,2,\cdots,
\end{eqnarray}
with $\mu_n$ given by (\ref{nwet}).

Note that the orthogonality measure is positive definite if $M<0$ (i.e. the value $M_0$ of the concentrated mass is positive).

In the following we determine the three term recurrence relation
that the polynomials $\tilde{P}_n^{(-1)}(x)$ verify. It is already
known that the little $q$-Jacobi polynomials satisfy the relation
\begin{eqnarray}
P_{n+1}+b_nP_n+u_nP_{n-1}=xP_{n}, \label{Jtr1}
\end{eqnarray}
and that the recurrence coefficients are defined by
\[
u_n=A_{n-1}C_n,\quad b_n=A_n+C_n,
\]where $A_n,C_n$ are given as
\begin{align}
A_n=q^n\frac{(1-aq^{n+1})(1-abq^{n+1})}{(1-abq^{2n+1})(1-abq^{2n+2})},\quad
C_n=aq^n\frac{(1-q^n)(1-bq^{n})}{(1-abq^{2n+1})(1-abq^{2n})}.
\end{align}
Under the Geronimus transformation
\begin{eqnarray}
\tilde{P}_{n}(x)={P}_n(x)-B_n{P}_{n-1}(x),
\end{eqnarray}
with $B_n=\frac{\Phi_n}{\Phi_{n-1}}$ and $\Phi_n$ defined by
\eqref{phiexp}, one obtains the three term recurrence relation
\begin{eqnarray}
\tilde{P}_{n+1}+\tilde{b}_n\tilde{P}_n+\tilde{u}_n\tilde{P}_{n-1}=x\tilde{P}_{n},\label{tp}
\end{eqnarray}
with coefficients
\begin{eqnarray}
&&\tilde{u}_1=\frac{\phi_1}{\phi_0^2}\qquad
\tilde{u}_n=\frac{u_{n-1}B_n}{B_{n-1}},\quad n=2,3,\cdots\\
&& \tilde{b}_0=b_0+\frac{\phi_1}{\phi_0}\qquad
\tilde{b}_n=b_n+B_{n+1}-B_n,\quad n=1,2,\cdots
\end{eqnarray}

When we set $a=q^2, q=-e^{\epsilon},b=-e^{\beta\epsilon}$ and take
the limit $\epsilon\rightarrow 0$, Eq. (\ref{tp}) reduces to
\begin{eqnarray}
\tilde
{P}_{n+1}^{(-1)}+\tilde{b}_n^{(-1)}\tilde{P}_n^{(-1)}+\tilde{u}_n^{(-1)}\tilde{P}_{n-1}^{(-1)}
=x\tilde{P}_{n}^{(-1)}.\label{recur}
\end{eqnarray}
Here the coefficients are
\begin{eqnarray}
&& \tilde{b}_n^{(-1)}=\lim\limits_{\epsilon\rightarrow
0}(b_n+B_{n+1}-B_n)=b_n^{(-1)}+\lim\limits_{\epsilon\rightarrow
0}(B_{n+1}-B_n) \\
&& \tilde{u}_n^{(-1)}=\lim\limits_{\epsilon\rightarrow
0}\frac{u_{n-1}B_n}{B_{n-1}}=u_{n-1}^{(-1)}\lim\limits_{\epsilon\rightarrow
0}\frac{B_n}{B_{n-1}},
\end{eqnarray}
where
\begin{eqnarray}
u_n^{(-1)}=-\frac{n(n+2)}{(2n+1+\beta)(2n+3+\beta)},\quad
b_n^{(-1)}=1
\end{eqnarray}
when $n$ is even, and
\begin{eqnarray}
u_n^{(-1)}=-\frac{(n+\beta)(n+2+\beta)}{(2n+1+\beta)(2n+3+\beta)},\quad
b_n^{(-1)}=-1
\end{eqnarray}
when $n$ is odd. Note that
\begin{align}
\lim\limits_{\epsilon\rightarrow 0}B_n=\lim\limits_{q\rightarrow
-1}\frac{\Phi_n}{\Phi_{n-1}}=\left\{\begin{matrix}
\frac{n+2}{2n+1+\beta}\frac{M-[(3+\beta)(1+\beta)]/[(n+2)(n+1+\beta)]}
{M-[(3+\beta)(1+\beta)]/[n(n+1+\beta)]}\quad  \texttt{n even}\\
\\
-\frac{n+2+\beta}{2n+1+\beta}
\frac{M-[(3+\beta)(1+\beta)]/[(n+1)(n+2+\beta)]}{M-[(3+\beta)(1+\beta)]/[(n+1)(n+\beta)]}
\quad \texttt{n odd}\quad.
\end{matrix}\right.
\end{align}

As a final observation, let us identify the matrix orthogonal
polynomials that the even (or odd) part of the
$\tilde{P}_n^{(-1)}(x)$ define. Split the polynomials
$\tilde{P}_n^{(-1)}(x)$ into its even ($E_n$) and odd ($O_n$) parts:
\begin{eqnarray}
\tilde{P}_n^{(-1)}(x)=E_n(x)+O_n(x).
\end{eqnarray}
From the recurrence relation (\ref{recur}), we have
\begin{align}
x^2 E_n=E_{n+2}&+(\tilde{b}_{n+1}^{(-1)}+
\tilde{b}_n^{(-1)})E_{n+1}+(\tilde{u}_{n+1}^{(-1)}
+\tilde{u}^{(-1)}_n+(\tilde{b}^{(-1)}_n)^2)E_n\nonumber\\
&
+(\tilde{b}^{(-1)}_n+\tilde{b}^{(-1)}_{n-1})\tilde{u}^{(-1)}_nE_{n-1}
+\tilde{u}^{(-1)}_{n-1}\tilde{u}^{(-1)}_nE_{n-2}.
\end{align}
With the redefinition $E_n=\sigma_n F_n, \,
\sigma_n=\sqrt{\tilde{u}^{(-1)}_1\tilde{u}_2^{(-1)}\cdots
\tilde{u}_{n}^{(-1)}}$, it is easy to see that the polynomial $F_n$
satisfies the five-term recurrence relation,
\begin{eqnarray}
x^2
F_n(x)=c_{n,0}F_n+c_{n,1}F_{n-1}+c_{n+1,1}F_{n+1}+c_{n,2}F_{n-2}+c_{n+2,2}F_{n+2}
\end{eqnarray}
where the coefficients are
\begin{eqnarray}
c_{n,0}=(\tilde{u}_{n+1}^{(-1)}
+\tilde{u}^{(-1)}_n+(\tilde{b}^{(-1)}_n)^2),
c_{n,1}=(\tilde{b}_{n-1}^{(-1)}+\tilde{b}_n^{(-1)})\sqrt{\tilde{u}_n^{(-1)}},
c_{n,2}=\sqrt{\tilde{u}_n^{(-1)}\tilde{u}_{n-1}^{(-1)}}.
\end{eqnarray}
From the theorem in \cite{duran}, the matrix polynomials
$\{\textsf{P}_n(x)\}$ defined by
\begin{eqnarray}
\textsf{P}_n(x)=\left(
\begin{matrix}
R_{2,0}(F_{2n})(x) &  R_{2,1}(F_{2n})(x)\\
R_{2,0}(F_{2n+1})(x) & R_{2,1}(F_{2n+1})(x)
\end{matrix}
\right)
\end{eqnarray}
satisfy the matrix three term recurrence relation
\begin{eqnarray}
x\textsf{P}_n(x)=\textsf{D}_{n+1}\textsf{P}_{n+1}(x)
+\textsf{E}_n\textsf{P}_n(x)+\textsf{D}_n^{\ast}\textsf{P}_{n-1}(x)
\end{eqnarray}
where
\begin{eqnarray}
\textsf{D}_n=\left(
\begin{matrix}
c_{2n,2}&  0 \\
c_{2n,1} & c_{2n+1,2}
\end{matrix}
\right),\qquad \textsf{E}_n=\left(
\begin{matrix}
c_{2n,0}&  c_{2n+1,1} \\
c_{2n+1,1} &c_{2n+1,0}
\end{matrix}
\right).
\end{eqnarray}
The polynomials $R_{N,m}(p)(x)$ are defined by
\begin{eqnarray}
R_{N,m}(p)(x)=\sum\limits_{n=0}^K\frac{p^{(nN+m)}(0)}{(nN+m)!}x^n,
\end{eqnarray}
where $K$ is the integer part of $(\deg(p(x)) -m)/N$

\section{Conclusion}
To sum up, we have added to the exploration of $-1$ polynomials by
introducing some $-1$ Krall-Jacobi polynomials. We focused on the
simplest positive definite case. They have been obtained from
generalized $q-$Jacobi polynomials through a limiting procedure and
have remarkable features. Noteworthy is the fact that they obey a
third-order differential-difference eigenvalue equation involving
the reflection operator. The 3-term recurrence relation has also
been determined. Finally let us stress that this presents an
interesting example of OPs whose measure involves a discrete mass at
the center of the orthogonality interval (rather than at its
boundary which is more common). Note also that the -1 Krall polynomials introduced 
in this paper have also appeared in a recent study of 
generalized Hahn-Geronimus classical orthogonal polynomials \cite{VZ_HG}

\section*{\bf Acknowledgements}
A.Zh. wishes to thank the Centre de Recherches Math\'{e}matiques
(CRM) at the Universit\'{e} de Montr\'{e}al for its hospitality. The
work of L.V. is supported by a grant from the National Science and
Engineering Research Council (NSERC) of Canada. G.F.Yu acknowledges
a postdoctoral fellowship from the Mathematical Physics Laboratory
of the CRM. He is also supported by the NNSF of China (Grant
no.11371251) and Chenguang Program (09CG08) sponsored by Shanghai
Municipal Education Commission and Shanghai Educational Development
Foundation. \vskip .5cm

\end{document}